\documentclass[a4paper,12pt]{article}
\usepackage{amsthm,amsmath}
\usepackage{amssymb}
\usepackage[all]{xypic}
\usepackage{enumerate}
\usepackage{a4wide}
\usepackage{hyperref}
\usepackage{latexsym}

\newtheorem{theorem}{Theorem}[section]
\newtheorem{proposition}[theorem]{Proposition}
\newtheorem{corollary}[theorem]{Corollary}
\newtheorem{lemma}[theorem]{Lemma}
\newtheorem*{theorem*}{Theorem}
\newtheorem*{proposition*}{Proposition}
\newtheorem*{corollary*}{Corollary}
\newtheorem*{lemma*}{Lemma}
\theoremstyle{definition}
\newtheorem{definition}[theorem]{Definition}

\newtheorem{example}[theorem]{Example}
\newtheorem{remark}[theorem]{Remark}
\newtheorem*{remark*}{Remark}

\newtheorem*{definition*}{Definition}

\newcommand{\coring}[1]{\mathfrak{#1}}
\newcommand{\tensor}[1]{\otimes_{#1}}

\newcommand{\rcomod}[1]{\mathcal{M}^{#1}}
\newcommand{\rmod}[1]{\mathcal{M}_{#1}}
\newcommand{\lmod}[1]{{}_{#1}\mathcal{M}}

\newcommand{\lcomod}[1]{{}^{#1}\mathcal{M}}

\renewcommand{\hom}[3]{\mathrm{Hom}_{#1}(#2,#3)}

\newcommand{\rend}[2]{\mathrm{End}({#2}_{#1})}
\newcommand{\lend}[2]{\mathrm{End}({}_{#1}#2)}

\newcommand{\rcomatrix}[2]{#2^* \tensor{#1} #2}

\newcommand{\Nat}[1]{\mathcal{N}at(#1,#1)}

\newcommand{\dostensor}[3]{#1 \tensor{#2} #3}
\newcommand{\trestensor}[5]{#1 \tensor{#2} #3 \tensor{#4} #5}
\newcommand{\fourtensor}[7]{#1 \tensor{#2} #3 \tensor{#4} #5 \tensor{#6} #7}

\renewcommand{\Nat}{\mathbb{N}}

\newcommand{\Real}{\mathbb{R}}
\newcommand{\Co}{\mathbb{C}}
\newcommand{\cuat}{\mathbb{H}}

\newcommand{\h}{\hspace{-3pt}-\hspace{-3pt}}

\newcommand{\otimescal}{\otimes_{\scriptscriptstyle{C(\alpha)}}}
\newcommand{\otimescbe}{\otimes_{\scriptscriptstyle{C(\beta)}}}
\newcommand{\otimesss}{\otimes_{\scriptscriptstyle{A_{\omega}(\alpha,\beta)}}}
\newcommand{\otimesmnk}{\otimes_{\scriptscriptstyle{M_n(k)}}}

\newcommand{\corb}{\Sigma^* \otimes_B \Sigma}
\newcommand{\corc}{\Sigma^* \otimes_C \Sigma}
\newcommand{\corge}{\mathfrak{C}}

\begin{document}
\title{Galois Corings and a Jacobson-Bourbaki type Correspondence.
}
\author{J. Cuadra\footnote{Supported by the grant BFM2002-02717 from MEC and FEDER.} \\
\normalsize Departamento de \'{A}lgebra y An\'{a}lisis Matem\'{a}tico \\
\normalsize Universidad de Almer\'{\i}a \\
\normalsize E04120 Almer\'{\i}a, Spain \\
\normalsize e-mail: \textsf{jcdiaz@ual.es} \and J. G\'omez-Torrecillas\footnote{Supported by the grant MTM2004-01406 from MEC and FEDER.} \\
\normalsize Departamento de \'{A}lgebra \\ \normalsize Universidad
de Granada\\ \normalsize E18071 Granada, Spain \\  \normalsize
e-mail: \textsf{gomezj@ugr.es} }

\date{}
\maketitle

\begin{abstract}
The Jacobson-Bourbaki Theorem for division rings was formulated
in terms of corings by Sweedler in \cite{Sweedler:1975}.
Finiteness conditions hypotheses are not required in this new
approach. In this paper we extend Sweedler's result to simple
artinian rings using a particular class of corings, comatrix
corings. A Jacobson-Bourbaki like correspondence for simple
artinian rings is then obtained by duality.
\end{abstract}

\section*{Introduction}

One of the key pieces in the Galois theory of fields and more
generally of division rings is the Jacobson-Bourbaki Theorem, see
\cite[Chapter 7, Sections 2, 3]{J2} and \cite[Section 8.2]{J}. Let
$E$ be a division ring with prime field $k$. Consider the
injective ring homomorphism $r:E \rightarrow \lend{k}{E}, e
\mapsto r_e$ where $r_e(e')=e'e$ for all $e' \in E$ (the
multiplication in $\lend{k}{E}$ is the opposite of the
composition). The Jacobson-Bourbaki Theorem states that there is a
bijective correspondence between the set of division subrings $D$
of $E$ such that $_DE$ is finite dimensional and the set of
subrings $S$ of $\lend{k}{E}$ such that $Im(r) \subseteq S$ and
$S_E$ is finite dimensional. The ring $\lend{k}{E}$ is indeed an
$E$-ring and the condition $Im(r) \subseteq S$ can be rephrased as
$S$ being an $E$-subring of $\lend{k}{E}$. This correspondence is
hidden behind the veil of the Galois connection in a Galois
extension of fields or more generally of division rings, see
\cite{J2} and \cite{Dr1,Dr2}.
\par \smallskip

Using the dual structure of $E$-ring, the structure of
$E$-coring, in \cite{Sweedler:1975} Sweedler gave a dual result
to the Jacobson-Bourbaki Theorem. The advantage of using this dual
structure is that the finiteness conditions needed in the
Jacobson-Bourbaki Theorem can be dropped. The finiteness
conditions come implicit in the structure of coring through the
fact that an element is mapped into a finite sum of elements via
the comultiplication map. Sweedler's result asserts that there is
a bijective correspondence between division subrings of $E$ and
quotient corings of the $E$-coring $E \otimes_k E$. The
Jacobson-Bourbaki Theorem can be obtained from Sweedler's result
by duality and this process makes clear why the finiteness
conditions are needed. \par \smallskip

The goal of this paper is to extend Sweedler's result from
division rings to simple artinian rings replacing the Sweedler
coring $E \otimes_k E$ by a more general type of coring, the
comatrix coring introduced in \cite{ElKaoutit/Gomez:2003} to
describe the structure of cosemisimple corings. Let $\Sigma$ be a
finitely generated and projective right module over a ring $A$ and
let $B$ be a simple artinian subring of $End(\Sigma_A)$. Let
$\coring{C}$ denotes the comatrix $A$-coring $\Sigma^* \otimes_B
\Sigma$ constructed on the bimodule $\Sigma$ and with coefficients
$B$, see (\ref{comco}). Our main theorem (Theorem \ref{galarti})
states that there is a bijective correspondence between the set of
all simple artinian subrings $B \subseteq C \subseteq
End(\Sigma_A)$ and the set of all coideals $J$ of $\coring{C}$
such that the quotient coring $\coring{C}/J$ is simple
cosemisimple. If in addition $\Sigma_A$ is simple, then any
quotient coring of $\coring{C}$ is simple cosemisimple, thus
obtaining a bijective correspondence between intermediate division
subrings of $B \subseteq \rend{A}{\Sigma}$ and coideals of
$\coring{C}$ (see Remark \ref{sigmasimple}). Sweedler's result,
together with some additional information on conjugated
subextensions, is then obtained as a consequence by taking $A$ a
division ring and $\Sigma=A$ (Corollary \ref{PJB}). An example
illustrating the bijective correspondence is worked out. This
closes Section \ref{predual}, which thus contains our main
results. \par \smallskip

Section \ref{comatrix} recalls the most fundamental results on
comatrix corings, Galois corings, and cosemisimple corings needed
in the sequel. We also include a homological characterization of
Galois corings (Theorem \ref{Galois}), which gives as a
consequence that if the canonical map is surjective for a
quasi-projective comodule with a generating condition, then the
coring is Galois (Corollary \ref{surinj}). This corollary is used
in the proof of the main result of Section \ref{dual}, namely,
Theorem \ref{JB}, which states a Jacobson-Bourbaki correspondence
for simple artinian subextensions of a ring extension. This
correspondence is dual to the stated in Section \ref{predual}. We
complete the paper with an Appendix that contains a complete
classification of the simple cosemisimple $\Co /
\mathbb{R}$--corings. \par \bigskip

We next fix notation and present some basic definitions. In the
sequel $A,B$ denote associative and unitary algebras over a
commutative ring $K$. By $\otimes_A$ we denote the tensor product
over $A$. The category of right $A$-modules is denoted by ${\cal
M}_A$. Bimodules are assumed to be centralized by $K$. An
$A$-coring (or $A/K$--coring, when $K$ is not obvious from the
context) is a triple $(\coring{C},\Delta,\epsilon)$ where
$\coring{C}$ is an $A$-bimodule and $\Delta:\coring{C} \rightarrow
\coring{C} \otimes_A\coring{C}$ (comultiplication) and
$\epsilon:\coring{C} \rightarrow A$ (counit) are $A$-bimodule maps
such that $(id_{\coring{C}} \otimes_A \Delta) \Delta=(\Delta
\otimes_A id_{\coring{C}})\Delta$ and $(\epsilon \otimes_A
id_{\coring{C}})\Delta=(id_{\coring{C}} \otimes_A
\epsilon)\Delta=id_{\coring{C}}$. In a categorical language, a
coring is just a coalgebra in the monoidal category of
$A$-bimodules with the tensor product $\otimes_A$ as a product.
For $c \in \coring{C}$ we will write $\Delta(c)=\sum c_{(1)}
\otimes_A c_{(2)}$. The left dual
${^*}\coring{C}=Hom({_A}\coring{C},{_A}A)$ of the coring
$\coring{C}$ is an $A^{opp}$-ring ($A^{opp}$ denotes the opposite
ring of $A$) with the product $(f*g)(c)=\sum f(c_{(1)}g(c_{(2)}))$
for all $f,g \in {^*}\coring{C}$ and $c \in \coring{C}$.
Similarly, the right dual $\coring{C}^*$ of $\coring{C}$ is an
$A^{opp}$-ring with the product $(f*g)(c) = \sum
g(f(c_{(1)})c_{(2)})$ for $f, g \in \coring{C}^*$ and $c \in
\coring{C}$. \par \smallskip

A right $\coring{C}$-comodule is a right $A$-module together with
an $A$-module map $\rho_M:M \rightarrow M \otimes_A \coring{C}$
such that $(id_M \otimes_A \Delta)\rho_M=(\rho_M \otimes
id_{\coring{C}})\rho_M$ and $(id_M \otimes \epsilon)\rho_M=id_M$.
A $\coring{C}$-comodule map between two right
$\coring{C}$-comodules $M$ and $N$ is an $A$-module map $f:M
\rightarrow N$ such that $(f \otimes_A
id_{\coring{C}})\rho_M=\rho_N f$. By $Hom_{\coring{C}}(M,N)$ we
will denote the $K$-module of all $\coring{C}$-comodule maps
between $M$ and $N$. The category whose objects are right
$\coring{C}$-comodules and whose morphisms are
$\coring{C}$-comodule maps is denoted by ${\cal M}^{\coring{C}}$.
It is an additive $K$--linear category and if $_A\coring{C}$ is
flat, then it is a Grothendieck category. The product of every
endomorphism ring of an object in an additive category is by
default the composition. We adopt, however, the following
convention in the case of modules: the product of the endomorphism
ring $\rend{A}{M}$ of a right $A$--module $M$ is the composition,
although by $\lend{A}{N}$ we will denote the opposite ring of the
endomorphism ring of a left $A$--module $N$, being then its
product the opposite of the composition.

\section{Comatrix corings,  Galois comodules and cosemisimple corings}\label{comatrix}

Let $\Sigma$ be a $B-A$--bimodule, and assume that $\Sigma_A$ is
finitely generated and projective with a finite dual basis
$\{(e_i^*,e_i) \} \subseteq \Sigma^* \times \Sigma$. We can
consider a coring structure \cite[Proposition
2.1]{ElKaoutit/Gomez:2003} over the $A$--bimodule
$\rcomatrix{B}{\Sigma}$ with comultiplication and counit defined
respectively by
\begin{equation}\label{comco}
\Delta(\dostensor{\phi}{B}{x}) = \sum_i
\fourtensor{\phi}{B}{e_i}{A}{e_i^*}{B}{x}, \qquad
\epsilon(\dostensor{\phi}{B}{x}) = \phi(x).
\end{equation}
The comultiplication is independent of the choice of the dual
basis. This coring will be called the $A$-comatrix coring on
$\Sigma$ with coefficients in $B$. The $A$-module $\Sigma$ becomes
a right $\rcomatrix{B}{\Sigma}$--comodule with coaction
\[
\varrho_{\Sigma}(x) = \sum_i\trestensor{e_i}{A}{e_i^*}{B}{x}.
\]

Assume $\Sigma$ to be the underlying $A$--module of a right
comodule over some $A$--coring $\coring{C}$, with structure map
$\rho_{\Sigma} : \Sigma \rightarrow
\dostensor{\Sigma}{A}{\coring{C}}$. In this case, with $T =
\rend{\coring{C}}{\Sigma}$, we have from \cite[Proposition
2.7]{ElKaoutit/Gomez:2003} that the map $\mathsf{can} :
\rcomatrix{T}{\Sigma} \rightarrow \coring{C}$ defined by
\[
\mathsf{can}(\dostensor{\phi}{T}{x}) = \sum \phi(x_0)x_1 \qquad
(\rho_{\Sigma}(x) = \sum \dostensor{x_0}{A}{x_1})
\]
is a homomorphism of $A$--corings. This canonical map allowed to
extend \cite[Definition 3.4]{ElKaoutit/Gomez:2003} the notion of a
Galois coring without assuming the existence of group-like
elements. When the role of the comodule $\Sigma$ is stressed, the
terminology of Galois comodules introduced in
\cite{Brzezinski/Wisbauer:2003} is more convenient. Probably, the
best solution here is to mention both the coring and the
comodule.

\begin{definition}
The pair $(\coring{C},\Sigma)$ is said to be \emph{Galois} if
$\mathsf{can}$ is an isomorphism. In such a case, we say that
$\coring{C}$ is a \emph{Galois coring} and $\Sigma$ is termed a
\emph{Galois comodule}. The extension $T \subseteq
\rend{\coring{C}}{\Sigma}$ is called a \emph{
$(\coring{C},\Sigma)$-Galois extension.}
\end{definition}

The notion of a noncommutative $G$--Galois extension, may be
recovered from this definition, see \cite[Example
2.9]{ElKaoutit/Gomez:2003}. In this case, the corresponding Galois
coring has a group-like element. We give an example of a
noncommutative Galois extension for a coring without group-like
elements.

\begin{example}{\emph Let $\Co$ and $\cuat$ denote the complex number field and
the Hamilton's quaternions algebra respectively. Consider the
right $\Co$-vector space $\coring{T}$ with basis $\{c,s\}$. This
vector space becomes a $\Co$-bimodule with left action $ic=c$ and
$is=-s$ and a $\Co$-coring with comultiplication and counity
defined by
$$\begin{array}{ll}
\Delta(c)=c \otimes c-s \otimes s, & \epsilon(c)=1, \vspace{3pt} \\
\Delta(s)=c \otimes s+s \otimes c, & \epsilon(s)=0.
\end{array}$$
Analogously to the coalgebra case this coring is called {\it the
trigonometric coring}. This coring has no group-like elements. Let
$\Sigma$ be a right $\Co$-vector space with basis $\{v_1,v_2\}$.
The map $\rho:\Sigma \rightarrow \Sigma \otimes \coring{T}$
defined by $$v_1 \mapsto v_1 \otimes c+v_2 \otimes s, \ v_2
\mapsto v_2 \otimes c-v_1 \otimes s$$ makes $\Sigma$ into a right
$\coring{T}$-comodule. It is not difficult to check that
$(\coring{T},\Sigma)$ is a Galois $\Co$-coring. The corresponding
Galois extension is the well-known embedding of $\cuat$ into
$M_2(\Co)$:
$$i \mapsto \left(\begin{array}{cc} 0 & -1 \\ 1 & 0 \end{array}\right),
\qquad j \mapsto \left(\begin{array}{cc} i & 0 \\ 0 & -i
\end{array}\right).$$
}
\end{example}
\vspace{0.6cm}

Our first objective is to enrich
\cite[18.26]{Brzezinski/Wisbauer:2003} with a new characterization
of Galois comodules. We need two previous observations. The first
one is that if $\coring{C}$ is flat as a left $A$--module then,
using that the forgetful functor $U : \rcomod{\coring{C}}
\rightarrow \rmod{A}$ is faithful and exact, the following lemma
can be proved (see \cite{Gomez/Zarouali:2004}).

\begin{lemma}\label{fg}
If $\coring{C}$ is flat as a left $A$--module, then a right
$\coring{C}$--comodule $M$ is finitely generated in the
Grothendieck category $\rcomod{\coring{C}}$ if and only if $M$ is
finitely generated as a right $A$--module.
\end{lemma}

 We have a pair of functors
\begin{equation}\label{adj1}
 \xymatrix{\rmod{T} \ar@<0.5ex>^{-\tensor{T}\Sigma}[rr] & &
\rcomod{\coring{C}} \ar@<0.5ex>^{\hom{\coring{C}}{\Sigma}{-}}[ll]}
\end{equation}
where $-\tensor{T} \Sigma$ is left adjoint to
$\hom{\coring{C}}{\Sigma}{-}$. If $\chi :
\dostensor{\hom{\coring{C}}{\Sigma}{-}}{T}{\Sigma} \rightarrow
id_{\rcomod{\coring{C}}}$ is the counit of this adjunction, then
the canonical map can be expressed \cite[Lemma
3.1]{ElKaoutit/Gomez:2003} as the composite
\begin{equation}\label{canchi}
\xymatrix{\mathsf{can} : \rcomatrix{T}{\Sigma}
\ar^-{\dostensor{(\dostensor{-}{A}{id_{\coring{C}}})\circ
\rho_{\Sigma}}{T}{\Sigma}}[rrr] & & &
\dostensor{\hom{\coring{C}}{\Sigma}{\coring{C}}}{T}{\Sigma}
\ar^-{\chi_{\coring{C}}}[r] & \coring{C}}
\end{equation}
Observe that $\chi_{\coring{C}}$ is an isomorphism if and only if
$\mathsf{can}$ is so. Our second observation is the following
lemma.

\begin{lemma}\label{etacan}
Let $\eta : id_{\rmod{T}} \rightarrow
\hom{\coring{C}}{\Sigma}{\dostensor{-}{T}{\Sigma}}$ be the unit of
the adjunction \eqref{adj1}. Consider the map
$\mathsf{can}_*:\hom{\coring{C}}{\Sigma}{\rcomatrix{T}{\Sigma}}
\rightarrow \hom{\coring{C}}{\Sigma}{\coring{C}}, \ f \mapsto
\mathsf{can} \circ f.$ Then the following composition is the
identity map on $\Sigma^*$
\[
\xymatrix{\Sigma^* \ar^-{\eta_{\Sigma^*}}[rr] & &
\hom{\coring{C}}{\Sigma}{\rcomatrix{T}{\Sigma}}
\ar^-{\mathsf{can}_*}[rr] & & \hom{\coring{C}}{\Sigma}{\coring{C}}
\cong \Sigma^*}
\]
\end{lemma}
\begin{proof}
If we apply the displayed composite map to $\phi \in \Sigma^*$,
then we obtain the map from $\Sigma$ to $A$ given by $x \mapsto
\sum \epsilon_{\coring{C}}(\phi(x_0)x_1)$, which is nothing but
$\phi$ since
\[
\sum \epsilon_{\coring{C}}(\phi(x_0)x_1) = \sum
\phi(x_0)\epsilon_{\coring{C}}(x_1) = \sum \phi (x_0
\epsilon_{\coring{C}}(x_1)) = \phi (x)
\]
\end{proof}

We are now in position to state our homological characterization
of Galois comodules.

\begin{theorem}\label{Galois}
Let $\Sigma$ be a right comodule over an $A$--coring $\coring{C}$
and assume that $\Sigma$ is finitely generated and projective as
a right $A$-module. If ${}_A\coring{C}$ is flat, then
$(\coring{C}, \Sigma)$ is Galois if and only if there exists an
exact sequence $\Sigma^{(J)} \rightarrow \Sigma^{(I)} \rightarrow
\coring{C} \rightarrow 0$ in $\rcomod{\coring{C}}$ such that the
sequence
\[
\xymatrix{\hom{\coring{C}}{\Sigma}{\Sigma^{(J)}} \ar[r] &
\hom{\coring{C}}{\Sigma}{\Sigma^{(I)}} \ar[r] &
\hom{\coring{C}}{\Sigma}{\coring{C}} \ar[r] & 0}
\]
is exact.
\end{theorem}
\begin{proof}
Consider $T^{(J)} \rightarrow T^{(I)} \rightarrow \Sigma^*
\rightarrow 0$ a free presentation of the left $T$--module
$\Sigma^*$. By tensorizing with ${}_T\Sigma$ we obtain an exact
sequence $\Sigma^{(J)} \rightarrow \Sigma^{(I)} \rightarrow
\rcomatrix{T}{\Sigma} \rightarrow 0.$ We have then the following
commutative diagram in $\rmod{T}$
\begin{equation}\label{etacaza}
\xymatrix{T^{(J)} \ar[r] \ar^-{\eta_{T^{(J)}}}[d] &  T^{(I)}
\ar[r] \ar^-{\eta_{T^{(I)}}}[d]& \Sigma^* \ar[r]
\ar^-{\eta_{\Sigma^*}}[d]
&  0 \\
\hom{\coring{C}}{\Sigma}{\Sigma^{(J)}} \ar[r]  &
\hom{\coring{C}}{\Sigma}{\Sigma^{(I)}} \ar[r] &
\hom{\coring{C}}{\Sigma}{\rcomatrix{T}{\Sigma}}  \ar[r] & 0}
\end{equation}
Now, $\eta_{T^{(J)}}$ and $\eta_{T^{(I)}}$ are isomorphisms
because $\Sigma$ is finitely generated in the Grothendieck
category $\rcomod{\coring{C}}$ (see \cite{Gomez/Zarouali:2004}).
By Lemma \ref{etacan}, $\eta_{\Sigma^*}$ is an isomorphism if and
only if $\mathsf{can}_*$ is bijective. So if we assume that
$(\coring{C},\Sigma)$ is Galois, then $\eta_{\Sigma^*}$ is an
isomorphism and, from \eqref{etacaza}, the following sequence is
exact
\[
\xymatrix{\hom{\coring{C}}{\Sigma}{\Sigma^{(J)}} \ar[r]  &
\hom{\coring{C}}{\Sigma}{\Sigma^{(I)}} \ar[r] &
\hom{\coring{C}}{\Sigma}{\rcomatrix{T}{\Sigma}}  \ar[r] & 0}
\]
Observe that the isomorphism of $A$--corings $\mathsf{can} :
\rcomatrix{T}{\Sigma} \cong \coring{C}$ is also an isomorphism of
right $\coring{C}$--comodules. This finishes the proof of the
necessary condition. For the sufficiency, consider the commutative
diagram in $\rcomod{\coring{C}}$ with exact rows
\[
\xymatrix{\Sigma^{(J)} \ar[r] & \Sigma^{(I)} \ar[r] & \coring{C}
\ar[r] & 0 \\
\dostensor{\hom{\coring{C}}{\Sigma}{\Sigma^{(J)}}}{T}{\Sigma}
\ar[r] \ar^-{\chi_{\Sigma^{(J)}}}[u]&
\dostensor{\hom{\coring{C}}{\Sigma}{\Sigma^{(I)}}}{T}{\Sigma}
\ar[r] \ar^-{\chi_{\Sigma^{(I)}}}[u]&
\dostensor{\hom{\coring{C}}{\Sigma}{\coring{C}}}{T}{\Sigma} \ar[r]
\ar^-{\chi_{\coring{C}}}[u] & 0}
\]
In this diagram, $\chi_{\Sigma^{(J)}}$ and $\chi_{\Sigma^{(I)}}$
are isomorphisms because $\Sigma$ is finitely generated in
$\rcomod{\coring{C}}$. We have then that $\chi_{\coring{C}}$ is an
isomorphism and, thus, $(\coring{C},\Sigma)$ is Galois.
\end{proof}

T. Brzezi\'{n}ski has shown in \cite{Brzezinski:2004} that a
simple comodule with $\mathsf{can}$ surjective is Galois. As a
consequence of Theorem \ref{Galois} we derive a generalization of
Brzezi\'{n}ski's result. Following the definition given for
modules in \cite{Anderson/Fuller:1992}, we say the comodule
$\Sigma$ is \emph{quasi-projective} if for every exact sequence
$\Sigma \rightarrow N \rightarrow 0$ in $\rcomod{\coring{C}}$ then
the sequence of abelian groups $\hom{\coring{C}}{\Sigma}{\Sigma}
\rightarrow \hom{\coring{C}}{\Sigma}{N} \rightarrow 0$ is exact.
Since, by Lemma \ref{fg}, $\Sigma$ is finitely generated in
$\rcomod{\coring{C}}$, a straightforward adaptation of
\cite[Proposition 16.2.(2)]{Anderson/Fuller:1992} to Grothendieck
categories gives that $\hom{\coring{C}}{\Sigma}{-}$ will already
preserve exact sequences of the form $\Sigma^{(I)} \rightarrow N
\rightarrow 0$. The following corollary is then easily deduced
from Theorem \ref{Galois}, making use once more of the
``\emph{AB5}'' condition.

\begin{corollary}\label{surinj}
Assume that $\Sigma$ is quasi-projective in $\rcomod{\coring{C}}$
and generates every subcomodule of any finite direct sum of copies
of $\Sigma$ (e.g., $\Sigma$ is a semisimple comodule). If
$\mathsf{can}$ is surjective, then $(\coring{C},\Sigma)$ is
Galois.
\end{corollary}

Now, let us recall from \cite{ElKaoutit/Gomez/Lobillo:2002,
ElKaoutit/Gomez:2003} the structure of cosemisimple corings, which
is tightly related to the coring version of the Generalized
Descent Theorem formulated in \cite[Theorem
3.10]{ElKaoutit/Gomez:2003}. A coring is said to be
\emph{cosemisimple} if it satisfies one of the equivalent
conditions of the following theorem.

\begin{theorem}\cite[Theorem 3.1]{ElKaoutit/Gomez/Lobillo:2002}
The following assertions for an $A$-coring $\mathfrak{C}$
\vspace{-1ex} are equivalent:
\begin{enumerate}
\itemsep 0pt \item[(i)] Every left $\coring{C}$-comodule is
semisimple and $\lcomod{\coring{C}}$ is abelian. \item[(ii)] Every
right $\coring{C}$-comodule is semisimple and
$\rcomod{\coring{C}}$ is abelian.
\item[(iii)]${_\coring{C}}\coring{C}$ is semisimple and
$\coring{C}_A$ is flat. \item[(iv)] $\coring{C}_{\coring{C}}$ is
semisimple and $_A\coring{C}$ is flat.
\end{enumerate}
\end{theorem}

A coring is called \emph{simple} if it has no non trivial
subbicomodules. It was proved in \cite[Theorem
3.7]{ElKaoutit/Gomez/Lobillo:2002} that any cosemisimple coring
decomposes in a unique way as a direct sum of simple cosemisimple
corings. An example of simple cosemisimple coring is the comatrix
$A$-coring $\Sigma^* \otimes_B \Sigma$, where $\Sigma_A$ is
finitely generated and projective and $B\subseteq End(\Sigma_A)$
is simple artinian. The following result shows that indeed all
simple cosemisimple corings can be obtained in this way.

\begin{proposition}\cite[Proposition 4.2]{ElKaoutit/Gomez:2003}
Let $\mathfrak{C}$ be a simple cosemisimple $A$-coring and
$\Sigma_\mathfrak{C}$ a finitely generated right
$\mathfrak{C}$-comodule. Then $T=End(\Sigma_\mathfrak{C})$ is
simple artinian, $\Sigma_A$ is finitely generated projective and
the canonical map $\mathsf{can}:\Sigma^* \otimes_T \Sigma
\rightarrow \mathfrak{C}$ is an isomorphism.
\end{proposition}

A more precise description of simple cosemisimple corings is given
by the following structure theorem. It may be viewed as a
generalization of the Artin-Wedderburn Theorem.

\begin{theorem}\cite[Theorem 4.3]{ElKaoutit/Gomez:2003}
An $A$-coring $\mathfrak{C}$ is simple cosemisimple if and only
if there is a finitely generated projective right $A$-module
$\Sigma$ and a division subring $D \subseteq End(\Sigma_A)$ such
that $\mathfrak{C} \cong \Sigma^* \otimes_D \Sigma$ as
$A$-corings. \par \smallskip

In such a case, if $\Gamma$ is another finitely generated and
projective right $A$-module and $E \subseteq End(\Gamma_A)$ is a
division subring, then $\mathfrak{C} \cong \Gamma^* \otimes_E
\Gamma$ if and only if there is an isomorphism of right
$A$-modules $g:\Sigma \rightarrow \Gamma$ such that $gDg^{-1}=E.$
\end{theorem}

In view of this structure theorem, for a field extension $A/k$,
the classification of simple cosemisimple $A$-corings centralized
by $k$ is reduced to the classification of finite dimensional
division algebras over $k$ and to the study of how these division
algebras embed in matrix algebras over $A$. The first problem
leads to the Brauer group theory of a field and the second one
can be treated with the help of the Skolem-Noether Theorem. The
complete classification for the field extension $\Co/\Real$ is
obtained in the Appendix.

\section{The Galois connection from a coring point of
view}\label{predual}

Let $\Sigma_A$ be a finitely generated and projective right
$A$-module and denote by $S = \rend{A}{\Sigma}$ its endomorphism
ring. Let $B \subseteq S$ be a subring, and consider the comatrix
$A$--coring $\coring{C} = \rcomatrix{B}{\Sigma}$. A typical
situation is to consider $B = k$, a field, and $S = M_n(k)$, the
ring of square matrices of order $n$ over $k$. Let
$\mathsf{Subext}(S/B)$ denote the set of all ring subextensions $B
\subseteq C \subseteq S$, and $\mathsf{Coideals}(\coring{C})$ be
the set of all coideals of $\coring{C}$. Consider the maps
\begin{equation}\label{galoisconnection}
\xymatrix{\mathcal{J} : \mathsf{Subext}(S/B) \ar@/^/[r] &
\mathsf{Coideals}(\coring{C}) : \mathcal{R} \ar@/^/[l]}
\end{equation}
defined as follows. For each subextension $C \in
\mathsf{Subext}(S/B)$ we have a canonical homomorphism of
$A$--corings $\coring{C} = \rcomatrix{B}{\Sigma} \rightarrow
\rcomatrix{C}{\Sigma}$, whose kernel $\mathcal{J}(C)$ is a coideal
of $\coring{C}$. Conversely, given a coideal $J$ of $\coring{C}$,
then, by \cite[Proposition 2.5]{ElKaoutit/Gomez:2003}, $B
\subseteq \rend{\coring{C}}{\Sigma} \subseteq
\mathrm{End}(\Sigma_{\corge/J}) \subseteq S$, so that
$\mathcal{R}(J) = \rend{\coring{C}/J}{\Sigma}$ is a subextension
of $B \subseteq S$. Both $\mathcal{J}$ and $\mathcal{R}$ are
inclusion preserving maps. The following proposition, which
generalizes \cite[Proposition 6.1]{Schauenburg:1996}, collects
some more of their relevant general properties.

\begin{proposition}\label{galois}
The maps defined in \eqref{galoisconnection} enjoy the following
properties.\vspace{-1ex}
\begin{enumerate}[(1)]
\itemsep 0pt \item \label{RJ} $\mathcal{R}\mathcal{J}(C) \supseteq
C$ for every $C \in \mathsf{Subext}(S/B)$. \item \label{JR}
$\mathcal{J}\mathcal{R}(J) \subseteq J$ for every $J \in
\mathsf{Coideals}(\coring{C})$. \item \label{RJI}
$\mathcal{R}\mathcal{J}(C) = C$ if and only if
$\rend{\rcomatrix{C}{\Sigma}}{\Sigma} = C$. \item \label{JRI}
$\mathcal{J}\mathcal{R}(J) = J$ if and only if
$(\coring{C}/J,\Sigma)$ is Galois. \item \label{GalRJ} The maps
$\mathcal{J}$ and $\mathcal{R}$ establish a bijection between the
set consisting of ring subextensions $B \subseteq C \subseteq S$
such that $\rend{\rcomatrix{C}{\Sigma}}{\Sigma} = C$ and the set
of coideals $J$ of $\coring{C}$ such that $(\coring{C}/J,\Sigma)$
is Galois.
\end{enumerate}
\end{proposition}

\begin{proof} A direct computation gives that $\mathcal{R}\mathcal{J}(C) =
\rend{\rcomatrix{C}{\Sigma}}{\Sigma}$ for every $C \in
\mathsf{Subext}(S/B)$. Thus,  \emph{(\ref{RJ})} follows from
\cite[Proposition 2.5]{ElKaoutit/Gomez:2003}. This gives also
\emph{(\ref{RJI})}. Now, for a given coideal $J \in
\mathsf{Coideals}(\coring{C})$, we have the commutative diagram
with exact rows
\begin{equation}\label{diagal}
\xymatrix{0 \ar[r] & \mathcal{JR}(J) \ar[d] \ar[r] &
\rcomatrix{B}{\Sigma} \ar@{=}[d] \ar[r] &
\rcomatrix{\rend{\coring{C}/J}{\Sigma}}{\Sigma}
\ar^{\mathsf{can}_{\coring{C}/J}}[d] \ar[r] & 0 \\
0 \ar[r] & J \ar[r] & \coring{C} \ar[r] & \coring{C}/J \ar[r] & 0}
\end{equation}
which implies \emph{(\ref{JR})} and \emph{(\ref{JRI})}. Finally,
let us prove \emph{(\ref{GalRJ})}: for any ring subextension $C
\in \mathsf{Subext}(S/B)$, the $A$-coring $\corc$ is Galois by
\cite[Lemma 3.9]{ElKaoutit/Gomez:2003}. Hence
$\corge/\mathcal{J}(C)$ is Galois. Conversely, for any $J\in
\mathsf{Coideals}(\coring{C})$ we have
$\mathrm{End}(\Sigma_{\corge}) \subseteq
\mathrm{End}(\Sigma_{\corge/J})$. If $\corge/J$ is Galois, then
the canonical map $\mathsf{can}_{\coring{C}/J}$ is an isomorphism
and, therefore,
$\mathrm{End}(\Sigma_{\Sigma^*\otimes_{\mathrm{End}(\Sigma_{\corge/J})}
\Sigma})=\mathrm{End}(\Sigma_{\corge/J})$. Statement
\emph{(\ref{GalRJ})} follows now from \emph{(\ref{RJI})} and
\emph{(\ref{JRI})}.
\end{proof}

\begin{remark}{\normalfont
If $_A(\corc)$ is locally projective (see e.g.
\cite[42.10]{Brzezinski/Wisbauer:2003} for this notion), then
$C=\mathrm{End}(\Sigma_{\corc})$ if and only if ${_C}\Sigma$ is
faithfully balanced, i.e., $C$ is isomorphic to the biendomophisms
ring of ${}_C\Sigma$ under the natural map. This is because for
$_A(\corc)$ locally projective,
$\mathrm{End}(\Sigma_{\corc})=\mathrm{End}(_{^*(\corc)}\Sigma)$,
see \cite[19.3]{Brzezinski/Wisbauer:2003}. Then, by
\cite[Proposition 2.1]{ElKaoutit/Gomez:2003},
${}^*(\rcomatrix{C}{\Sigma}) \cong \lend{C}{\Sigma}^{op}$
canonically and, therefore,
$$\begin{array}{ll}
\mathrm{End}(\Sigma_{\corc}) &
=\mathrm{End}(_{^*(\corc)}\Sigma)=\mathrm{End}({_{\mathrm{End}({_C}\Sigma)^{op}}}\Sigma)
= \mathrm{End}(\Sigma_{\mathrm{End}({_C}\Sigma)}).\end{array}$$ }
\end{remark}

\begin{remark}
We have from \cite[Proposition 2.5]{ElKaoutit/Gomez:2003} that
\[
\rend{\rcomatrix{C}{\Sigma}}{\Sigma} = \{ f \in \rend{A}{\Sigma}
\; | \; \dostensor{f}{C}{x} = \dostensor{1}{C}{f(x)}, \hbox{ for
every } x \in \Sigma \} := \overline{C}.
\]
Since $\rcomatrix{C}{\Sigma}$ is Galois, we get that
$\overline{\overline{C}} = \overline{C}$. Proposition \ref{galois}
gives a bijective correspondence between coideals $J$ of
$\coring{C}$ such that $\coring{C}/J$ is Galois and subextensions
$C$ such that $C = \overline{C}$.
\end{remark}

We are now ready to state a generalization of Sweedler's predual
to Jacobson-Bourbaki Theorem \cite[Theorem 2.1]{Sweedler:1975}. We
will say that $C, C' \in \mathsf{Subext}(S/B)$ are
\emph{conjugated} in $S$ if there is an unit $g \in S$ such that
$C' = gCg^{-1}$.

\begin{theorem}\label{galarti}
Let $\Sigma$ be a finitely generated projective right $A$-module
and let $S=\mathrm{End}(\Sigma_A)$. Let $B$ be a subring of $S$
and consider the comatrix $A$-coring $\mathfrak{C}=\Sigma^*
\otimes_B \Sigma$. Denote by $\mathcal{S}$ the set of all simple
artinian subrings $B \subseteq C \subseteq S$ and let
$\mathcal{T}$ denote the set of all coideals $J$ of $\mathfrak{C}$
such that $\mathfrak{C}/J$ is simple cosemisimple. Then the maps
$$\begin{array}{l}
\mathcal{R}(-):\mathcal{T} \rightarrow \mathcal{S}, \ J \mapsto
\mathrm{End}(\Sigma_{\corge/J}) \vspace{3pt} \\
\mathcal{J}(-):\mathcal{S} \rightarrow \mathcal{T},\ C \mapsto
Ker(\corge \twoheadrightarrow \corc)
\end{array}$$
are inverse to each other. If, in addition, $A$ is a
(noncommutative) local ring, then two intermediate simple artinian
subrings $C$ and $C'$ are conjugated in $S$ if and only if
$\corge/\mathcal{J}(C)$ and $\corge/\mathcal{J}(C')$ are
isomorphic as $A$--corings.
\end{theorem}

\proof For $C \in \mathcal{S}$, the $A$-coring $\Sigma^* \otimes_C
\Sigma$ is a simple cosemisimple $A$-coring in virtue of
\cite[Proposition 4.2]{ElKaoutit/Gomez:2003}. Hence
$\corge/\mathcal{J}(C) \cong \corc$ is simple cosemisimple and so
$\mathcal{J}(C) \in \mathcal{T}$. From \cite[Proposition
2.5]{ElKaoutit/Gomez:2003}, $C \subseteq
\mathrm{End}(\Sigma_{\corc}).$ Since $_C\Sigma$ is faithfully
flat, \cite[Theorem 3.10]{ElKaoutit/Gomez:2003} yields
$C=\mathrm{End}(\Sigma_{\corc})$. By Theorem \ref{galois},
$\mathcal{R}\mathcal{J}(C) = C$.
\par \smallskip

Assume that $J \in \mathcal{T}$, then $\corge/J$ is simple
cosemisimple. By \cite[Theorem 4.1]{ElKaoutit/Gomez:2003},
$\corge/J$ is flat as a left $A$-module which implies, by
\cite[Lemma 3.1]{Gomez/Zarouali:2004}, that $\Sigma$ is finitely
generated as a right $\corge/J$-comodule. Hence
$\mathrm{End}(\Sigma_{\corge/J})$ is a simple artinian ring. So
$\mathcal{R}(J) \in \mathcal{S}$. Furthermore, $\corge/J$ is
Galois by \cite[Proposition 4.2]{ElKaoutit/Gomez:2003}. By Theorem
\ref{galois}, $\mathcal{J}\mathcal{R}(J) = J$. \par
\smallskip

For the second assertion, let $C,C'$ be intermediate simple
artinian subrings and let $g \in S$ be invertible such that
$gCg^{-1}=C'$. We check that $\Sigma^* \otimes_{C} \Sigma$ is
isomorphic to $\Sigma^* \otimes_{C'} \Sigma$. The set
$\{e_i^*g,g(e_i)\}_{i=1}^n$ is a dual basis for $\Sigma$ and, by
\cite[Remark 2.2]{ElKaoutit/Gomez:2003} the coring structure on
$\Sigma^* \otimes_{C'} \Sigma$ defined by this dual basis is the
same as the one defined by $\{e_i^*,e_i\}_{i=1}^n.$ Consider the
$A$-bimodule map,
$$\psi:\Sigma^* \times \Sigma \rightarrow \Sigma^* \otimes_{C'}
\Sigma, \ (\varphi,x) \mapsto \varphi g^{-1} \otimes_{C'} g(x).$$
Let $c \in C$ and $c'=gcg^{-1} \in C'$. Then,
$$\begin{array}{ll}
\psi(\varphi c,x) & = \varphi c g^{-1} \otimes_{C'} g(x) = \varphi
g^{-1}c'\otimes_{C'} g(x) = \varphi g^{-1} \otimes_{C'} c'g(x) \\
& =\varphi g^{-1} \otimes_{C'} gc(x)= \psi(\varphi,cx).
\end{array}$$
Hence $\psi$ defines a unique $A$-bimodule map
$\bar{\psi}:\Sigma^* \otimes_{C} \Sigma \rightarrow \Sigma^*
\otimes_{C'} \Sigma$. It is routine to verify that $\psi$ is an
isomorphism of $A$-corings. Conversely, let $\chi:\Sigma^*
\otimes_{C} \Sigma \rightarrow \Sigma^* \otimes_{C'} \Sigma$ be an
isomorphism of $A$-corings. Let $S,S'$ be the unique, up to
isomorphism, simple right $\Sigma^* \otimes_{C} \Sigma$-comodule
and $\Sigma^* \otimes_{C'} \Sigma$-comodule, respectively. Then
$\Sigma \cong S^{(m)}$ as a $\Sigma^* \otimes_{C} \Sigma$-comodule
and $\Sigma \cong S'^{(n)}$ as a $\Sigma^* \otimes_{C'}
\Sigma$-comodule for some $m,n \in \Nat$. Let $S^{\chi}$ denote
$S$ viewed as a right $\Sigma^* \otimes_{C'} \Sigma$-comodule via
$\chi$. Then $S^{\chi} \cong S'$ as a $\Sigma^* \otimes_{C'}
\Sigma$-comodule. In particular, they are isomorphic as right
$A$-modules. Since $A$ is local, $S^{\chi} \cong S' \cong A^{(l)}$
as right $A$-modules for a certain $l \in \Nat$. Let
$\Sigma^{\chi}$ denote $\Sigma$ when considered as a right
$\rcomatrix{C'}{\Sigma}$--comodule. Then $\Sigma^{\chi} \cong
(S^{\chi})^{(m)} \cong A^{(lm)}$ as a right $A$-module. On the
other hand, $\Sigma \cong A^{(ln)}$ as a right $A$-module. Since
the underlying right $A$-module of $\Sigma^{\chi}$ and $\Sigma$ is
the same, $m=n$ and hence $\Sigma^{\chi} \cong \Sigma$ as a right
$\Sigma^* \otimes_{C'} \Sigma$-comodule. Denote this isomorphism
by $g$. Then $\mathrm{End}(\Sigma_{\Sigma^* \otimes_{C'} \Sigma})=
\mathrm{End}(\Sigma^{\chi}_{\Sigma^* \otimes_{C'} \Sigma})$ via $d
\mapsto g^{-1}dg$. As $\mathrm{End}(\Sigma^{\chi}_{\Sigma^*
\otimes_{C'} \Sigma})=\mathrm{End}(\Sigma_{\Sigma^* \otimes_{C}
\Sigma}),$ $C=\mathrm{End}(\Sigma_{\Sigma^* \otimes_{C} \Sigma})$
and $C'=\mathrm{End}(\Sigma_{\Sigma^* \otimes_{C'} \Sigma}),$ the
assertion holds. \qed

\begin{remark}\label{sigmasimple}{\normalfont
With hypothesis as in Theorem \ref{galarti}, if we assume in
addition that $\Sigma_A$ is simple, then any quotient coring of
$\mathfrak{C}=\Sigma^* \otimes_B \Sigma$ is simple cosemisimple:
since $\Sigma_A$ is simple, $\Sigma$ is simple as a right
$\mathfrak{C}/J$-comodule for any coideal $J$ of $\mathfrak{C}$.
It is easy to see that the canonical map $\mathsf{can}:\Sigma^*
\otimes_{\mathrm{End}(\Sigma_{\mathfrak{C}/J})} \Sigma \rightarrow
\mathfrak{C}/J$ is surjective. Using that
$\Sigma_{\mathfrak{C}/J}$ is simple, by Corollary \ref{surinj},
$\mathsf{can}$ is an isomorphism. As
$\mathrm{End}(\Sigma_{\mathfrak{C}/J})$ is a division ring,
$\mathfrak{C}/J$ is simple cosemisimple. Thus we have a bijective
correspondence between intermediate division subrings of $B
\subseteq \mathrm{End}(\Sigma_A)$ and coideals of $\mathfrak{C}$.
Observe that no assumptions are made on $A$. As a consequence we
can derive Sweedler's predual to the Jacobson-Bourbaki Theorem,
with the additional information concerning conjugated subrings.}
\end{remark}

\begin{corollary}\label{PJB}\cite[Theorem 2.1]{Sweedler:1975}
Let $D \subseteq E$ be division rings. Set
$\mathfrak{C}=E\otimes_D E$ and let $g=1 \otimes_D 1$ be the
distinguished group-like element. For a coideal $J$ of
$\mathfrak{C}$ let $\pi_J:\mathfrak{C} \rightarrow \mathfrak{C}/J$
denote the canonical projection. Then, the maps
$$\begin{array}{l}
\mathcal{R}(-):\mathcal{T} \rightarrow \mathcal{S}, \ J \mapsto
\{e \in E: e\pi_J(g)=\pi_J(g)e\} \vspace{3pt} \\
\mathcal{J}(-):\mathcal{S} \rightarrow \mathcal{T},\ C \mapsto
Ker(\corge \twoheadrightarrow \corc)
\end{array}$$
establish a bijective correspondence between the set $\mathcal{S}$
of intermediate division rings $D \subseteq C \subseteq E$ and the
set $\mathcal{T}$ of coideals $J$ of $\mathfrak{C}$. Moreover, two
intermediate division rings $C$ and $C'$ are conjugated in $S$ if
and only if $\corge/\mathcal{J}(C) \cong \corge/\mathcal{J}(C').$
\end{corollary}

\proof It only remains to prove that
$\mathrm{End}(E_{\corge/J})=\{e \in E : e\pi_J(g)=\pi_J(g)e\}$ but
this is easily checked. \qed

\begin{remark}{\normalfont
If $\Sigma_A$ is not simple, then the quotient corings of $\corb$
need not in general to be simple cosemisimple. Let $k$ be a field,
$\Sigma=k^{(2)}$ and $T=T_2(k) \subseteq \lend{k}{\Sigma} =
M_2(k)$ the upper triangular matrix algebra. Then $\Sigma^*
\otimes_T \Sigma$ is a non simple cosemisimple quotient coalgebra
of $\Sigma^* \otimes_k \Sigma$. This example also serves to show
that factor corings of Galois corings are not Galois. The module
$\Sigma_T$ is isomorphic to the indecomposable projective $eT$,
where $e \in M_2(k)$ is the elementary matrix with $1$ in the
$(1,1)$-entry and zero elsewhere. Thus, $\rend{T}{\Sigma} \cong
eTe \cong k$, and the canonical map \eqref{canU} gives here a
surjective $k$-coalgebra homomorphism
$\mathsf{can}:\rcomatrix{k}{\Sigma} \rightarrow T^*$ which,
obviously, cannot be bijective. Thus, if we take $\coring{C} =
\rcomatrix{k}{\Sigma}$, and $\coring{D} = T^*$, then the factor
coalgebra $(\coring{D},\Sigma)$ of the Galois coalgebra
$(\coring{C},\Sigma)$ is not Galois. On the other hand, observe
that $\Sigma_{T^*}$ is projective but does not generate all its
submodules, which sheds some light on the conditions involved in
Theorem \ref{Galois}. }
\end{remark}

\begin{example}
We next illustrate the Galois connection established in
Theorem \ref{galarti} by a concrete example. Assume that $k$ has
an $n$-th primitive root of unity $\omega$. For $\alpha,\beta$ non
zero elements in $k$ let $A_{\omega}(\alpha,\beta)$ denote the
associative $k$-algebra generated by two elements $x,y$ subject to
the relations $x^n=\alpha,y^n=\beta$ and $yx=\omega xy$. Details
on the properties of this algebra to be used in the sequel may be
consulted in \cite[Chapter 15]{Mi}. The algebra
$A_{\omega}(\alpha,\beta)$ is a central simple $k$-algebra. For
our purposes we will assume that the subalgebras
$C(\alpha)=k\langle x:x^n=\alpha\rangle$ and $C(\beta)=k\langle
y:y^n=\beta\rangle$ are fields. Let $\Sigma$ be an $n$-dimensional
$C(\alpha)$-vector space with basis $B=\{v_1,...,v_n\}$ and
consider a dual basis $B^*=\{v_1^*,...,v_n^*\}$ in $\Sigma^*.$ The
algebra $A_{\omega}(\alpha,\beta)$ can be embedded in
$M_n(C(\alpha))$ by assigning
$$x\mapsto X=xe_{1,1}+\omega xe_{2,2}+...+\omega^{n-1}xe_{n,n},\qquad
y \mapsto Y=e_{1,2}+...+e_{n-1,n}+\beta e_{n,1},$$ where $e_{i,j}$
denotes the elementary matrix in $M_n(C(\alpha))$ with $1$ in the
$(i,j)$-entry and zero elsewhere. The action of $X$ and $Y$ on the
bases $B$ and $B^*$ is:
\begin{equation}\label{eqaction}
\begin{tabular}{lp{1cm}l}
$X\cdot v_j=\omega^{j-1}v_jx$ &  & $v_j^*\cdot
X=\omega^{j-1}xv_j^*$ \vspace{3pt} \\
$Y \cdot v_j=\left\{\begin{array}{ll} \beta v_{n} & \textrm{if $j=1$} \\
v_{j-1} & \textrm{if $j>1$}
\end{array}\right.$ & & $v_j^* \cdot Y=\left\{\begin{array}{ll}
\beta v_{1}^* & \textrm{if $j=n$} \\
v_{j+1}^* & \textrm{if $j<n$}
\end{array}\right.$
\end{tabular}
\end{equation}
If either $\alpha$ or $\beta$ is equal to $1$, then the algebra
$A_{\omega}(\alpha,\beta)$ is isomorphic to $M_n(k)$. \par
\medskip

We will next describe the coideals of the $C(\alpha)$--coring
$\Sigma^* \otimes_k \Sigma$ corresponding to the intermediate
extensions of $k \subset M_n(C(\alpha))$ given in the following
diagram:
$$
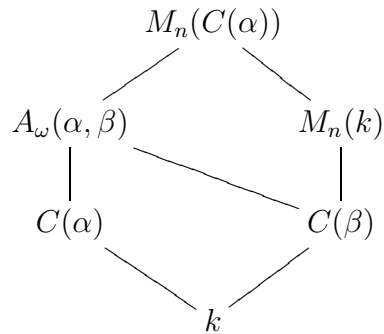
\xymatrix@R=0.7cm@C=0.0cm{ & M_n(C(\alpha)) &  \\
A_{\omega}(\alpha,\beta) \ar@{-}[ur] &  & M_n(k) \ar@{-}[ul] \\
C(\alpha) \ar@{-}[u]& & \ar@{-}[ull] C(\beta) \ar@{-}[u]  \\
 &  \ar@{-}[ul] k \ar@{-}[ur]  & }$$
For $l,m,i=1,...,n$ set $z_{l,m}^i=\alpha^{-1}x^{n-i}v_l^*
\otimes_k v_m x^i$. Observe that the set
$\{z_{l,m}^i:l,m,i=1,...,n\}$ is a basis of $\Sigma^* \otimes_k
\Sigma$ as a right $C(\alpha)$-vector space. We have that
$xz_{l,m}^i=z_{l,m}^{i-1}x$ for $1\leq i \leq n$ with the
convention $z_{l,m}^0=z_{l,m}^n.$ The comultiplication and counit
of $\Sigma^* \otimes_k \Sigma$ reads:
\begin{equation}\label{eqcomuco}
\begin{array}{ll}
\Delta(z_{l,m}^i) & = \sum_{j=1}^n (\alpha^{-1}x^{n-i}v_l^*
\otimes_k v_j) \otimescal (v_j^* \otimes_k v_m x^i) \vspace{3pt}
\\
 & = \sum_{j=1}^n (\alpha^{-1}x^{n-i}v_l^*
\otimes_k v_jx^{i}) \otimescal (\alpha^{-1}x^{n-i}v_j^* \otimes_k
v_m x^i) \vspace{3pt} \\
 & = \sum_{j=1}^n z_{l,j}^i \otimescal z_{j,m}^i, \vspace{3pt} \\
\epsilon(z_{l,m}^i) & =\delta_{l,m}.
\end{array}
\end{equation}

{\it Trivial extensions:} The trivial extensions $k \subseteq k
\subset M_n(C(\alpha))$ and $k \subset M_n(C(\alpha)) \subseteq
M_n(C(\alpha))$ correspond to the coideals $\{0\}$ and
$Ker(\epsilon)$ respectively. \par \medskip

{\it Extension $k \subset C(\alpha) \subset M_n(C(\alpha)):$} We
embed $C(\alpha)$ into $M_n(C(\alpha))$ by mapping $x$ to $X$ and
consider the $C(\alpha)$-coring $\Sigma^* \otimescal \Sigma$. The
action of $X$ on $B$ and $B^*$ gives the following relations in
$\Sigma^* \otimescal \Sigma:$
\begin{equation}\label{eqalpha}
\alpha^{-1}x^{n-i}v_l^*\otimescal v_m x^i=v_l^* \otimescal v_m.
\end{equation}
The set $\{v_l^* \otimescal v_m: l,m=1,...,n\}$ is a basis of
$\Sigma^* \otimescal \Sigma$ as a right $C(\alpha)$-vector space.
Set $c_{l,m}=v_l^* \otimescal v_m$. The bimodule structure on
this coring is given by $xc_{l,m}=\omega^{m-l}c_{l,m}x$. The
comultiplication and counit in $\Sigma^* \otimescal \Sigma$ are
defined by:
$$\Delta(c_{l,m})=\sum_{j=1}^n c_{l,j} \otimescal c_{j,m}, \qquad
\epsilon(c_{l,m})=\delta_{l,m}.$$ The coideal $J_{C(\alpha)}$ of
$\Sigma^* \otimes_k \Sigma$ corresponding to this extension is
the right subspace generated by the set
$$\{z_{l,m}^n-\omega^{i(m-l)}z_{l,m}^i:l,m,i=1,...,n\h1\}.$$
Observe that if we embed diagonally $C(\alpha)$ into
$M_n(C(\alpha))$, then $\Sigma$ is centralized by $C(\alpha)$ and
hence $\Sigma^* \otimescal \Sigma$ is indeed a coalgebra, the
comatrix coalgebra over $C(\alpha)$ of order $n$. Hence the
diagonal embedding of $C(\alpha)$ is not conjugated with the
preceding one. \par
\medskip

{\it Extension $k \subset C(\beta) \subset M_n(C(\alpha)):$} We
embed $C(\beta)$ into $M_n(C(\alpha))$ by mapping $y$ to $Y$ and
consider the $C(\alpha)$-coring $\Sigma^* \otimescbe \Sigma$.
Taking into account the action of $Y$ on $B$ and $B^*$, the
following relations in $\Sigma^* \otimescbe \Sigma$ are obtained:
\begin{equation}\label{eqbeta}
v_i^* \otimescbe v_j=\left\{\begin{array}{ll}
\beta^{-1}v^*_{n-(j-i)+1}\otimescbe v_1 & \textrm{if $i<j$} \\
v^*_{i-j+1} \otimescbe v_1 & \textrm{if $i \geq j$}
\end{array}\right.
\end{equation}
A basis of $\Sigma^* \otimescbe \Sigma$ as a left
$C(\alpha)$-vector space is:
$$\{\alpha^{-1}x^{n-i}v^*_l \otimescbe
v_1x^i:l,i=1,...,n\}.$$ Setting $c_{l}^i=\alpha^{-1}x^{n-i}v^*_l
\otimescbe v_1x^i$, the left action of $C(\alpha)$ on $\Sigma^*
\otimescbe \Sigma$ reads as $xc_l^i=c_l^{i-1}x$ with the
convention $c_l^0=c_l^{n}.$ The comultiplication and counit of
$\Sigma^* \otimescbe \Sigma$ is given by:
$$\begin{array}{ll}
\Delta(c_l^i) & = \sum_{j=1}^n (\alpha^{-1}x^{n-i}v^*_l \otimescbe
v_j) \otimescal (v_j^* \otimescbe v_1x^i) \vspace{3pt} \\
 & = \sum_{j=1}^n (\alpha^{-1}x^{n-i}v^*_l \otimescbe
v_jx^i) \otimescal (\alpha^{-1}x^{n-i}v_j^* \otimescbe v_1x^i) \vspace{3pt} \\
 & = \sum_{j=1}^l (\alpha^{-1}x^{n-i}v^*_l \otimescbe
v_jx^i) \otimescal (\alpha^{-1}x^{n-i}v_j^* \otimescbe v_1x^i) \vspace{3pt} \\
 & \quad + \sum_{j=l+1}^n (\alpha^{-1}x^{n-i}v^*_l \otimescbe
v_jx^i) \otimescal (\alpha^{-1}x^{n-i}v_j^* \otimescbe v_1x^i)
\vspace{3pt} \\
 & = \sum_{j=1}^l (\alpha^{-1}x^{n-i}v^*_{l-j+1} \otimescbe
v_1x^i) \otimescal (\alpha^{-1}x^{n-i}v_j^* \otimescbe v_1x^i) \vspace{3pt} \\
 & \quad + \beta^{-1}\sum_{j=l+1}^n (\alpha^{-1}x^{n-i}v^*_{n-(j-l)+1} \otimescbe
v_1x^i) \otimescal (\alpha^{-1}x^{n-i}v_j^* \otimescbe v_1x^i)
\vspace{3pt} \\
 & =\sum_{j=1}^l  c_j^i \otimescal c_{l-j+1}^i+
 \beta^{-1}\sum_{j=l+1}^n c_j^i \otimescal c_{n-(j-l)+1}^i , \vspace{5pt} \\

\epsilon(c_l^i) & =\delta_{l,1}.
\end{array}$$

The coideal $J_{C(\beta)}$ of $\Sigma^* \otimes_k \Sigma$
corresponding to $C(\beta)$ is the right subspace generated by the
following set:
$$\begin{array}{l}
\{z_{l,m}^i-\beta^{-1}z_{n-(m-l)+1,1}^i:l,m,i=1,...,n; l<m\} \cup
\{z_{l,m}^i-z_{l-m+1,1}^i:l,m,i=1,...,n; l \geq m\}.
\end{array}$$
\par \medskip

{\it Extension $k \subset A_{\omega}(\alpha,\beta) \subset
M_n(C(\alpha)):$} Consider $A_{\omega}(\alpha,\beta)$ as embedded
into $M_n(C(\alpha))$ by mapping $x$ to $X$ and $y$ to $Y$. We
next describe the $C(\alpha)$-coring $\Sigma^* \otimesss \Sigma$.
Since $C(\alpha)$ and $C(\beta)$ are contained in
$A_{\omega}(\alpha,\beta)$, similar relations to (\ref{eqalpha})
and (\ref{eqbeta}) are obtained. The set $\{v_i \otimesss
v_1:i=1,...,n\}$ is a basis of $\Sigma^* \otimesss \Sigma$ as a
right $C(\alpha)$-vector space. Set $c_i=v_i \otimesss v_1$. The
bimodule structure of $\Sigma^* \otimesss \Sigma$ is
$xc_i=\omega^{-i+1}c_ix$. The comultiplication and the counit of
$\Sigma^* \otimesss \Sigma$ are:
$$\begin{array}{ll}
\Delta(c_i) & = \sum_{l=1}^n (v_i^* \otimesss v_l)
\otimescal (v_l^* \otimesss v_1) \vspace{3pt} \\
 & = \sum_{l=1}^i (v^*_{i-l+1} \otimesss v_1) \otimescal (v_l^*
 \otimesss v_1) \vspace{3pt} \\
 & \quad + \beta^{-1}\sum_{l=i+1}^n (v^*_{n-(l-i)+1} \otimesss v_1) \otimescal (v_l^*
 \otimesss v_1), \vspace{3pt} \\
& = \sum_{l=1}^i c_l \otimescal c_{i-l+1} +
\beta^{-1}\sum_{l=i+1}^n c_{l}
\otimescal c_{n-(l-i)+1}, \vspace{5pt} \\
\epsilon(c_i) & =\delta_{i,1}.
\end{array}$$
The coideal $J_{A_{\omega}(\alpha,\beta)}$ of $\Sigma^* \otimes_k
\Sigma$ that corresponds to this intermediate extension is the
right subspace generated by the set:
$$\begin{array}{l}
\{z_{l,m}^n-\omega^{i(m-l)}z_{l,m}^i:l,m,i=1,...,n\}
\cup  \{z_{l,m}^i-\beta^{-1}z_{n-(m-l)+1,1}^i:l,m,i=1,...,n; \\
l<m\} \cup \{z_{l,m}^i-z_{l-m+1,1}^i:l,m,i=1,...,n; l \geq m\}.
\end{array}$$
\par \medskip

{\it Extension $k \subset M_n(k) \subset M_n(C(\alpha)):$} Let
$$X=e_{1,1}+\omega e_{2,2}+...+\omega^{n-1}e_{n,n},\qquad
Y=e_{1,2}+...+e_{n-1,n}+e_{n,1}.$$ Then $X^n=1,Y^n=1$ and
$YX=\omega XY$ and the $k$-algebra generated by $X$ and $Y$ is
$M_n(k)$. The action of $X$ and $Y$ on the basis $B$ and $B^*$ is
obtained from (\ref{eqaction}) for $\beta=1$. From these actions
we get the relations:
$$\begin{array}{ll}
v_1^* \otimesmnk v_1 & = v_n^* \cdot Y \otimesmnk v_1 = v_n^*
\otimesmnk Y \cdot v_1 =v_n^* \otimesmnk v_n  \vspace{3pt} \\
 & = v_{n-1}^* \cdot Y \otimesmnk v_n = v_{n-1}^*
\otimesmnk Y \cdot v_n =v_{n-1}^* \otimesmnk v_{n-1}  \vspace{3pt} \\
 & =.... \vspace{3pt} \\
 & = v_2^* \otimesmnk v_2. \vspace{3pt} \\
v_i^* \otimesmnk v_j & = \omega^{-j+1}v_i^* \otimesmnk X \cdot
v_j = \omega^{-j+1}v_i^*\cdot X \otimesmnk v_j \vspace{3pt} \\
 & =\omega^{i-j}v_i^* \otimesmnk v_j.
\end{array}$$
Then $v_i^* \otimesmnk v_j=0$ for $i \neq j$. Set
$c_i=\alpha^{-1}x^{n-i} v_1^* \otimesmnk v_1 x^{n-i}$ for
$i=1,...,n$. Then the set $\{c_i:i=1,...,n\}$ is a basis of
$\Sigma^* \otimesmnk \Sigma$ as a right $C(\alpha)$-vector space.
The bimodule structure of this coring is given by $xc_i=c_{i-1}x $
with the convention $c_0=c_n$. The comultiplication and counit of
$\Sigma^* \otimesmnk \Sigma$ takes the form:
$$\begin{array}{ll}
\Delta(c_i) & = \sum_{j=1}^n (\alpha x^{n-i}v_1^* \otimesmnk v_j)
\otimescal (v_j^* \otimesmnk v_1)x^i \vspace{3pt} \\
 & =(\alpha x^{n-i}v_1^* \otimesmnk v_1)
\otimescal (v_1^* \otimesmnk v_1)x^i \vspace{3pt} \\
 & = (\alpha x^{n-i}v_1^* \otimesmnk v_1x^i)
\otimescal (\alpha^{-1}x^{n-i} v_1^* \otimesmnk v_1)x^i \vspace{3pt} \\
 & = c_i \otimescal c_i, \vspace{3pt} \\
\epsilon(c_i) & =1.
\end{array}$$
This coring can also be obtained as the Sweedler coring associated
to the extension $k \subset C(\alpha)$. The coideal $J_{M_n(k)}$
of $\Sigma^* \otimes_k \Sigma$ corresponding to this extension is
the right subspace spanned by the set
$$\{z_{l,m}^i:l,m,i=1,...,n;l \neq m\}\cup
\{z_{1,1}^i-z_{l,l}^i:l,i=1,...,n\}.$$
\end{example}

\section{Duality}\label{dual}

Let $f:\coring{C} \rightarrow \coring{D}, g:\coring{C} \rightarrow
\coring{E}$ be surjective homomorphisms of $A$--corings. Then $Ker
f = Ker g$ if and only if there exists an isomorphism of
$A$--corings $\coring{D} \cong \coring{E}$ making commute the
diagram
\begin{equation}\label{triangulo}
\xymatrix{ & \coring{C} \ar_{f}[dl] \ar^{g}[dr] & \\
\coring{D} \ar^{\simeq}[rr] &  &  \coring{E}}
\end{equation}
Thus, every coideal $J$ of $\coring{C}$ determines a class of
surjective homomorphisms $\coring{C} \rightarrow \coring{D}$
having $J$ as their common kernel or, alternatively, the morphisms
in each class are connected by commutative triangles as in
\eqref{triangulo}.

From a formal point of view, corings over $A$ are dual to
$A$--rings, being these last understood to be morphisms of rings
$A \rightarrow U$. The definition of a homomorphism of $A$--rings
is obvious, and we will conceive an \emph{$A$--subring} of a given
$A$--ring $A \rightarrow E$ as an isomorphism class of injective
homomorphisms of $A$--rings $U \rightarrow E$. Obviously, every
$A$--subring of $E$ may be represented by an inclusion $U
\subseteq E$. We can thus consider the set $\mathsf{Subrings}(E)$
of $A$--subrings of $E$.

One of the possible concrete dual correspondences from
$A$--corings to $A$--rings goes as follows: if $\coring{C}$ is an
$A$--coring, then ${}^*\epsilon_{\coring{C}} : A \rightarrow
{}^*\coring{C}^{op}$ is an $A$--ring (see \cite[Proposition
3.2]{Sweedler:1975}), and under this mapping, homomorphisms of
$A$--corings give homomorphisms of $A$--rings. In particular, if
$J \in \mathsf{Coideals}(\coring{C})$ and $\coring{C} \rightarrow
\coring{D} = \coring{C}/J$ is the corresponding canonical
projection, then we have an injective homomorphism of $A$--rings
${}^*\coring{D}^{op} \rightarrow {}^*\coring{C}^{op}$ (see again
\cite[Proposition 3.2]{Sweedler:1975}). If $\coring{C} =
\rcomatrix{B}{\Sigma}$ is a comatrix $A$--coring, then, by
\cite[Proposition 2.1]{ElKaoutit/Gomez:2003}, we have an
injective homomorphism of $A$--rings ${}^*\coring{D}^{op}
\rightarrow {}^*\coring{C}^{op} \cong \lend{B}{\Sigma}$. The
corresponding $A$--subring of $\lend{B}{\Sigma}$ will be denoted
by $\mathcal{R}'(J)$. Conversely, given an $A$--subring $U
\rightarrow \lend{B}{\Sigma}$, then $\rend{U}{\Sigma}$ is
independent on the representative $U$ of the $A$--subring. We
define the coideal $\mathcal{J}'(U)$ of $\rcomatrix{B}{\Sigma}$
as the kernel of the homomorphism of $A$--corings
$\rcomatrix{B}{\Sigma} \rightarrow
\rcomatrix{\rend{U}{\Sigma}}{\Sigma}$ induced by the ring
homomorphism $B \rightarrow \rend{U}{\Sigma}$.

These considerations lead to a Galois connection for $E =
\lend{B}{\Sigma}$:
\begin{equation}\label{galoisconnection2}
\xymatrix{\mathcal{J}' : \mathsf{Subrings}(E) \ar@/^/[r] &
\mathsf{Coideals}(\coring{C}) : \mathcal{R}' \ar@/^/[l]}
\end{equation}

The following proposition collects some of its relevant
properties. A right comodule $\Sigma$ over an $A$--coring
$\coring{D}$ is said to be \emph{loyal} if the canonical map
\[
\xymatrix{\rcomatrix{\rend{\coring{D}}{M}}{M} \ar[r] &
\rcomatrix{\lend{{}^*\coring{D}}{M}}{M}}
\]
 induced by the inclusion
$\rend{\coring{D}}{\Sigma} \subseteq
\lend{{}^*\coring{D}}{\Sigma}$ is a bijection. By
\cite[19.2,19.3]{Brzezinski/Wisbauer:2003}, if ${}_A\coring{D}$ is
locally projective, then every right $\coring{D}$--comodule is
loyal.

\begin{proposition}
The mappings defined in \eqref{galoisconnection2} enjoy the
following properties. \vspace{-1ex}
\begin{enumerate}[(1)]
\itemsep 0pt
\item\label{R'J'} $\mathcal{R}'\mathcal{J}'(U) =
\lend{\rend{U}{\Sigma}}{\Sigma}$ for every $U \in
\mathsf{Subrings}(E)$ and, thus, $U \subseteq
\mathcal{R}'\mathcal{J}'(U)$. \item\label{J'R'} $J \supseteq
\mathcal{J}'\mathcal{R}'(J)$ for every $J \in
\mathsf{Coideals}(\rcomatrix{B}{\Sigma})$ such that
$\Sigma_{\coring{C}/J}$ is loyal. \item\label{dualGalois} The maps
$\mathcal{J}'$ and $\mathcal{R}'$ establish a bijection between
the set of $A$--subrings $U$ of $\lend{B}{\Sigma}$ such that
$\Sigma_U$ is faithful and balanced, and the set of coideals $J$
of $\coring{C}$ such that $(\coring{C}/J,\Sigma)$ is Galois and
$\Sigma_{\coring{C}/J}$ is loyal.
\end{enumerate}
\end{proposition}
\begin{proof}\emph{(\ref{R'J'})} Given an $A$--subring $U \subseteq
\rend{B}{\Sigma}$, $\mathcal{J}'(U)$ is defined as the kernel of
the surjective homomorphism of $A$--corings $\rcomatrix{B}{\Sigma}
\rightarrow \rcomatrix{\rend{U}{\Sigma}}{\Sigma}$. The commutative
diagram of injective homomorphisms of $A$--rings
\[
\xymatrix{{}^*(\rcomatrix{\rend{U}{\Sigma}}{\Sigma})^{op} \ar[r]
\ar^{\simeq}[d] & {}^*(\rcomatrix{B}{\Sigma})^{op}
\ar^{\simeq}[d] \\
\lend{\rend{U}{\Sigma}}{\Sigma} \ar[r] & \lend{B}{\Sigma}}
\]
deduced from \cite[Proposition 2.1]{ElKaoutit/Gomez:2003}, shows
that $\mathcal{R}'\mathcal{J}'(U) =
\lend{\rend{U}{\Sigma}}{\Sigma}$.\par \smallskip

\emph{(\ref{J'R'})} With $p: \coring{C} \rightarrow \coring{D} =
\coring{C}/J$ the canonical projection, consider the commutative
diagram
\begin{equation}\label{triangulo2}
\xymatrix{\rcomatrix{B}{\Sigma} \ar^{p}[rr] \ar^{f}[dr] & & \coring{D} \\
 & \rcomatrix{\rend{\coring{D}}{\Sigma}}{\Sigma} \ar^{\mathsf{can}}[ur]
 &},
\end{equation}
where $f$ is induced by the morphism $B \rightarrow
\rend{\coring{D}}{\Sigma}$. Since $\Sigma_{\coring{D}}$ is loyal,
we get from the diagram that $J \supseteq \ker (f) =
\mathcal{J}'\mathcal{R}'(J)$.\par \smallskip

\emph{(\ref{dualGalois})} Let $U \subseteq E = \lend{B}{\Sigma}$
be an $A$--subring. By definition, $\mathcal{J}'(U)$ is such that
$\coring{C}/\mathcal{J}'(U) \cong
\rcomatrix{\rend{U}{\Sigma}}{\Sigma}$. By \cite[Lemma
3.9]{ElKaoutit/Gomez:2003},
$(\rcomatrix{\rend{U}{\Sigma}}{\Sigma}, \Sigma)$ is Galois. This
means that the canonical map
\[
\xymatrix{\mathsf{can} :
\rcomatrix{\rend{\rcomatrix{\rend{U}{\Sigma}}{\Sigma}}{\Sigma}}{\Sigma}
\ar[r] & \rcomatrix{\rend{U}{\Sigma}}{\Sigma}}
\]
is an isomorphism. But this map is nothing but the one induced by
the ring extension
\[
\rend{\rcomatrix{\rend{U}{\Sigma}}{\Sigma}}{\Sigma} \subseteq
\lend{{}^*(\rcomatrix{\rend{U}{\Sigma}}{\Sigma})}{\Sigma} =
\rend{\rend{\rend{U}{\Sigma}}{\Sigma}}{\Sigma} = \rend{U}{\Sigma},
\]
which proves that $\Sigma$ is a loyal right
$\rcomatrix{\rend{U}{\Sigma}}{\Sigma}$--comodule. Part
\emph{(\ref{R'J'})} gives obviously $\mathcal{R}'\mathcal{J}'(U) =
U$. Conversely, let $J$ be a coideal of $\coring{C}$ such that
$(\coring{C}/J,\Sigma)$ is Galois and $\Sigma_{\coring{C}/J}$ is
Galois. Put $\coring{D} = \coring{C}/J$. From the triangle
\eqref{triangulo2}, and the fact that $\Sigma_\coring{D}$ is
loyal, we compute the $A$--subring $\mathcal{R}'(J)$ of
$\lend{R}{\Sigma}$ as
\begin{equation}
{}^*(\coring{D})^{op} \cong
{}^*(\rcomatrix{\rend{\coring{D}}{\Sigma}}{\Sigma})^{op} \cong
{}^*(\rcomatrix{\rend{{}^*\coring{D}^{op}}{\Sigma}}{\Sigma})^{op}
\cong \lend{\rend{{}^*\coring{D}^{op}}{\Sigma}}{\Sigma}
\end{equation}
This means that $\mathcal{R}'(J)$, defined as
${}^*\coring{D}^{op}$, is such that $\Sigma_{\mathcal{R}'(J)}$ is
faithful and balanced. The diagram \eqref{triangulo2} gives in
addition that $\mathcal{J}'\mathcal{R}'(J) = J$.
\end{proof}

We are now in position to prove our version for simple artinian
rings of the Jacobson-Bourbaki theorem.

\begin{theorem}\label{JB}
Let $\Sigma$ be a finitely generated projective right module over
a Quasi-Frobenius ring $A$ and let $S=\mathrm{End}(\Sigma_A)$. Let
$B$ be a subring of $S$. Denote by $\mathcal{I}$ the set of
intermediate simple artinian subrings $B \subseteq C \subseteq S$
such that $_C\Sigma$ is finitely generated. Let $\mathcal{D}$
denote the set of all $A$-subrings $U$ of
$\mathrm{End}({_B}\Sigma)$ such that $U_A$ is finitely generated
and projective, and $\Sigma_U$ is semisimple and isotypic. The
maps
$$\begin{array}{l}
\mathcal{R}^*(-):\mathcal{D} \rightarrow \mathcal{I}, \ U \mapsto
\mathrm{End}(\Sigma_U), \vspace{3pt} \\
\mathcal{J}^*(-):\mathcal{I} \rightarrow \mathcal{D}, \ C \mapsto
\mathrm{End}({_C}\Sigma).
\end{array}$$
establish a bijective correspondence. This correspondence is dual
to that of Theorem \ref{galarti}. Moreover, if $A$ is in addition
local, then two intermediate simple artinian rings $C$ and $C'$
are conjugated if and only if $\lend{C}{A}$ and $\lend{C'}{A}$ are
isomorphic as $A$--rings.
\end{theorem}

\begin{proof} If $C$ is simple artinian and ${}_C\Sigma$ is finitely
generated, then $U = \lend{C}{\Sigma}$ is a simple artinian
$A$--ring. Thus, $\Sigma_U$ is semisimple isotypic. The comatrix
$A$--coring $\rcomatrix{C}{\Sigma}$ is then finitely generated and
projective as a left $A$--module. The ring isomorphism $U^{op}
\cong {}^*(\rcomatrix{C}{\Sigma})$ given in \cite[Proposition
2.1]{ElKaoutit/Gomez:2003} is an isomorphism of $A$--bimodules. In
particular, $U_A$ becomes a finitely generated module.  Obviously,
${}_C\Sigma$ is faithful and balanced. Conversely, assume that
$\Sigma_U$ is semisimple isotypic for an $A$--subring $U \subseteq
\lend{B}{\Sigma}$ such that $U_A$ is finitely generated. Then $C =
\rend{U}{\Sigma}$ is a simple artinian subring of $S$, because
$\Sigma_U$ is clearly finitely generated. To prove that $\Sigma_U$
is faithful and balanced, consider that the comatrix $A$--coring
structure on $\rcomatrix{U}{U}$ induces, via the isomorphism
$\rcomatrix{U}{U} \cong U^*$ an $A$--coring structure on this last
$A$--bimodule \cite[Example 2.4]{ElKaoutit/Gomez:2003}. If $\{
u_{\alpha}^*, u_{\alpha}\}_{\alpha=1}^n$ is a finite dual basis
for $U_A$, the comultiplication and counit are given explicitly by
\[
\Delta : U^* \rightarrow U^*\tensor{A}U^*, \quad \varphi \mapsto
\sum_{\alpha} \varphi u_{\alpha} \tensor{A} u_{\alpha}^*,
\]
\[
\epsilon : U^* \rightarrow A, \quad \varphi \mapsto \varphi(1)
\]
The canonical map $U \rightarrow {}^*(U^*)$ is then an
anti-isomorphism of rings. We may then identify the categories of
$\rcomod{U^*}$ and $\lmod{U^{op}} = \rmod{U}$ (see e.g.
\cite[19.6]{Brzezinski/Wisbauer:2003}). In fact, an explicit
isomorphism of categories is given as follows.  For each element
$m$ in a right $U$--module $M$, the equality $mu =\sum_{\alpha}
mu_{\alpha}u_{\alpha}^*(u)$ for $u \in U$ says that
$\{mu_{\alpha}, u_{\alpha}^* \}$ is a set of right rational
parameters in the sense of \cite{ElKaoutit/Gomez/Lobillo:2002}.
Thus \cite[Corollary 4.7]{ElKaoutit/Gomez/Lobillo:2002} gives the
isomorphism of categories $\rmod{U} = \rcomod{U^*}$, where the
right $U^*$--comodule structure on $M_U$ is given by $\rho_M(m) =
\sum_{\alpha} mu_{\alpha} \tensor{A} u_{\alpha}^*$. In
particular, the (finitely generated) simple isotypic right
$U$--module may be considered as a right $U^*$--comodule, and we
have a canonical map given by
\begin{equation}\label{canU}
\mathsf{can} : \rcomatrix{C}{\Sigma} \rightarrow U^*, \quad
\varphi \tensor{C} x \mapsto \sum_{\alpha}
\varphi(xu_{\alpha})u_{\alpha}^*
\end{equation}
Now, for every $0 \neq u \in U$, let $x \in \Sigma$ such that $xu
\neq 0$. Then, for $\varphi \in \Sigma^*$ with $\varphi(xu) \neq
0$, we have
\[
 \mathsf{can}(\varphi \tensor{C} x)(u) =\sum_{\alpha}
 \varphi(xu_{\alpha})u_{\alpha}^*(u) =\sum_{\alpha}
 \varphi(xu_{\alpha}u_{\alpha}^*(u)) =\varphi (xu) \neq 0,
\]
which implies, being $A$ Quasi-Frobenius and $U_A$ finitely
generated, that $\mathsf{can}$ is surjective. Finally, since
$\Sigma_{U^*}$ is semisimple, we get from Corollary \ref{surinj}
that $\mathsf{can}$ an isomorphism. Therefore, we have an
isomorphism
\[
\xymatrix{U^{op} \cong {}^*(U^*) \ar^{\mathsf{can}^*}[r] &
\lend{C}{\Sigma}^{op}}
\]
which turns out to be, by using \cite[equation
(3)]{ElKaoutit/Gomez:2003}, the canonical homomorphism $U
\rightarrow \lend{C}{\Sigma}$. In this way, $\Sigma_U$ is faithful
and balanced and ${}_C\Sigma$ must be finitely generated.
\end{proof}

\begin{remark}
If the base ring $A$ is simple artinian then Theorem \ref{JB}
takes a simpler form. Thus, for a simple artinian ring $B$, the
maps
$$\begin{array}{l}
\mathcal{R}^*(-):\mathcal{D} \rightarrow \mathcal{I}, \ U \mapsto
End(\Sigma_U), \vspace{3pt} \\
\mathcal{J}^*(-):\mathcal{I} \rightarrow \mathcal{D}, \ C \mapsto
End({_C}\Sigma).
\end{array}$$
establish a bijective correspondence between the set $\mathcal{D}$
of simple artinian $A$-subrings $U$ of $End({_B}\Sigma)$ such that
$U_A$ is finitely generated and the set $\mathcal{I}$ of
intermediate simple artinian subrings $B \subseteq C \subseteq S$
such that $_CS$ is finitely generated.
\end{remark}

If $\Sigma = A$, and $B \subseteq A$ is a subring, then for every
$A$--subring $U \subseteq \lend{B}{A}$, we have that
$\mathcal{R}^* (U) = \rend{A}{U} = \{ c \in A : u(ca) = cu(a), u
\in U, a \in A \}$. If, moreover, $A$ is a division ring, then we
deduce from Theorem \ref{JB} the following version of the
Jacobson-Bourbaki theorem (c.f. \cite[\S 4]{Sweedler:1975}).

\begin{corollary}
Let $B\subseteq A$ be division rings. The maps $\mathcal{R}^*$
and $\mathcal{J}^*$ establish a bijective correspondence between
the intermediate division rings $B \subseteq C \subseteq A$ such
that $_CA$ is finite dimensional and the $A$--subrings $U
\subseteq \lend{B}{A}$ such that $U_A$ is finite dimensional.
Moreover, two intermediate division subrings $C$ and $C'$ are
conjugated if and only if $\lend{C}{A}$ and $\lend{C'}{A}$ are
isomorphic $A$--rings.
\end{corollary}
\begin{proof}
The pertinent remark here is that $A_U$ is always a simple module
(since $A_A$ is simple).
\end{proof}

\section*{Appendix}

Throughout this section all corings are $\Co/\mathbb{R}$-corings
and all $\Co$-bimodules are assumed to be centralized by the field
of real numbers $\Real$.
\par
\smallskip

We know by the Structure Theorem of simple and cosemisimple
corings that any such a $\Co$-coring $\corge$ is of the form
$\Sigma^* \otimes_{D} \Sigma$ where $\Sigma$ is a finite
dimensional complex vector space and $D$ is a $\Real$-division
algebra embedded in $\mathrm{End}(\Sigma_{\Co}) \cong M_n(\Co)$
where $n=dim(\Sigma_{\Co}).$ Furthermore, two corings  $\Sigma^*
\otimes_{D} \Sigma$ and $\Sigma^* \otimes_{E} \Sigma$ are
isomorphic if and only if there is an invertible $u \in
\mathrm{End}(\Sigma_{\Co})$ such that $uEu^{-1}=D$, i.e., $E$ and
$D$ are conjugated in $\mathrm{End}(\Sigma_{\Co})$. \par
\smallskip

By Fr\"{o}benius Theorem, $D=\Real,\Co$ or $\cuat$. We study the
possible ways of embedding these division algebras in $M_n(\Co)$.
As a consequence of the proof of the Skolem-Noether Theorem
(\cite[Theorem 4.9]{J}) the non conjugated ways of embedding $D$
in $M_n(\Co)$ are in bijective correspondence with the simple left
$D \otimes_{\Real} M_n(\Co)$-modules. \par \medskip

{\bf 1. Case $D=\Real$.} Since $M_n(\Co)$ is an $\Real$-algebra,
the only way of embedding $\Real$ in $M_n(\Co)$ is the obvious
one. The comultiplication and counit of the coring $\Sigma^*
\otimes_{\Real} \Sigma$ is described in (\ref{eqcomuco}). If $n=
1$, then this coring is isomorphic to the coring
$[\mathbb{Z}/2]\mathbb{C}$ associated to the canonical
$\mathbb{Z}/2$--grading on $\mathbb{C}$. As a right
$\mathbb{C}$--vector space, $[\mathbb{Z}/2]\mathbb{C}$ is free
over the basis $\mathbb{Z}/2 = \{ [0], [1] \}$. Its left
$\mathbb{C}$-vector space structure is determined by the rules
$i[0] = [1]i$, $i[1] = [0]i$. In this coring, $[0]$ and $[1]$ are
group-like elements. For $n > 1$, the coring
$\rcomatrix{\mathbb{R}}{\Sigma}$ is isomorphic with the tensor
product coring (see \cite[Proposition 1.5]{Gomez/Louly:2003})
$[\mathbb{Z}/2]\mathbb{C} \tensor{\mathbb{R}} M^c(\mathbb{R},n)$,
where $M^c(\mathbb{R},n)$ is the comatrix $\mathbb{R}$--coalgebra.
An explicit isomorphism sends $z_{m,l}^0$ onto $[0]
\tensor{\mathbb{R}} x_{m,l}$ and $z_{m,l}^1$ onto $[1]
\tensor{\mathbb{R}} x_{m,l}$, where $x_{m,l}$ denotes the matrix
with $1$ in the component $(l,m)$ and $0$ elsewhere.
\par
\smallskip

{\bf 2. Case $D=\Co$.} Since $\Co \otimes_{\Real} M_n(\Co) \cong
M_n(\Co) \oplus M_n(\Co)$, there are two ways (for $n > 0$) of
embedding $\Co$ into $M_n(\Co)$. Let $e_{p,q}$ denote the
elementary matrix in $M_n(\Co)$ with $1$ in the $(p,q)$-entry and
zero elsewhere. The two non conjugate embeddings are represented
by the one sending $i$ to $i(\sum_{l=1}^n e_{l,l})$ and the one
sending $i$ to $\bar{i}=\sum_{l=1}^n e_{l,n-l+1}.$ In the first
case, $\Sigma$ is centralized by $\Co$ and therefore $\Sigma^*
\otimes_{\Co} \Sigma$ is a $\Co$-coalgebra. Thus $\Sigma^*
\otimes_{\Co} \Sigma$ is isomorphic to the comatrix coalgebra of
order $n$. Let us study the second case. The action of $\bar{i}$
on $\Sigma$ and $\Sigma^*$ is:
$$\bar{i}\cdot v_j=v_{n-j+1}i,\qquad v_j^*\cdot
\bar{i}=iv_{n-j+1}^*.$$ Then the bimodule structure on $\Sigma^*
\otimes_{\Co} \Sigma$ is given by:
$$\begin{array}{ll}
i(v_p^* \otimes_{\Co} v_q) & = iv_p^* \otimes_{\Co}
v_q=v^*_{n-p+1} \cdot \bar{i} \otimes_{\Co} v_q = v^*_{n-p+1}
\otimes_{\Co} \bar{i} \cdot v_q \\
 & = v^*_{n-p+1} \otimes_{\Co} v_{n-q+1}i = (v^*_{n-p+1}
 \otimes_{\Co} v_{n-q+1})i.
\end{array}$$
The comultiplication and counit of $\Sigma^* \otimes_{\Co} \Sigma$
is defined by:
$$\Delta(v^*_p \otimes v_q)=\sum_{l=1}^n (v_p^* \otimes v_l)
\otimes (v_l^* \otimes v_q),\qquad \epsilon(v_p^* \otimes
v_q)=\delta_{p,q}.$$ This coring can be described as the right
$\Co$-vector space $\mathfrak{C}=\oplus_{p,q=1}^n v_{p,q}\Co$ with
bimodule structure $iv_{p,q}=v_{n-p+1,n-q+1}i$ and with
comultiplication and counit given by:
$$\Delta(v_{p,q})=\sum_{l=1}^n v_{p,l} \otimes_{\Co} v_{l,q},
\qquad \epsilon(v_{p,q})=\delta_{p,q}.$$
\par\smallskip

{\bf 3. Case $D=\cuat$.} Since $\cuat$ is a central simple
$\Real$-algebra and $M_n(\Co)$ is simple, $\cuat \otimes_{\Real}
M_n(\Co)$ is simple. Hence all embeddings of $\cuat$ in $M_n(\Co)$
are conjugate. Let us observe that if $\cuat$ embeds in
$M_n(\Co)$, then $n$ is even. By the Double Centralizer Theorem we
would have $M_n(\Co) \cong \cuat \otimes_{\Real} C(\cuat)$, where
$C(\cuat)$ denotes the centralizer of $\cuat$ in $M_n(\Co)$.
Comparing real dimensions, $4$ divides $2n^2$. Let
$\bar{i},\bar{j}$ be the generators of $\cuat$. Consider the
following embedding of $\cuat$ in $M_2(\Co):$
$$\bar{i} \mapsto \left(\begin{array}{cc} 0 & -1 \\ 1 & 0 \end{array}\right),
\qquad \bar{j} \mapsto \left(\begin{array}{cc} i & 0 \\ 0 & -i
\end{array}\right).$$
Since $n$ is even, the above embedding gives an embedding of
$\cuat$ in $M_n(\Co)$ by placing each block repeatedly in the main
diagonal. Any other embedding of $\cuat$ in $M_n(\Co)$ is
conjugated to this one. \par \smallskip

We next describe the coring $\Sigma^* \otimes_{\cuat} \Sigma$. We
first assume that $n=2$. Observe that this case is a particular
case of our example in the second section by taking $k=\Real$ and
$\alpha=\beta=-1$. The bimodule structure on $\Sigma^*
\otimes_{\cuat} \Sigma$ is the following:
$$i(v_1^* \otimes v_1)= (v^*_1 \otimes v_1)i, \qquad
i(v_2^* \otimes v_1)= -(v^*_2 \otimes v_1)i.$$ The
comultiplication and counit of $\Sigma^* \otimes_{\cuat} \Sigma$
read as:
$$\begin{array}{ll}
\Delta(v_1^* \otimes v_1)= (v_1^* \otimes v_1)\otimes (v_1^*
\otimes v_1)-(v_2^* \otimes v_1) \otimes (v_2^* \otimes v_1), &
\qquad \epsilon(v_1^* \otimes v_1)=1, \vspace{3pt} \\
\Delta(v_2^* \otimes v_1)= (v_2^* \otimes v_1)\otimes (v_1^*
\otimes v_1)+(v_1^* \otimes v_1) \otimes (v_2^* \otimes v_1), &
\qquad \epsilon(v_2^* \otimes v_1)=0.
\end{array}$$
This coring is precisely the {\it trigonometric coring}. We now
discuss the general case $dim(\Sigma_{\Co})=n=2m$. Let us recall
that $\cuat$ is embedded in $M_n(\Co)$ in the following way:
$$\begin{array}{ll}
\bar{i} \mapsto \sum_{i=1}^m e_{2l-1,2l}-\sum_{l=1}^m e_{2l,2l-1},
\qquad \bar{j} \mapsto i(\sum_{l=1}^n (-1)^{l+1}e_{l,l}).
\end{array}$$
Through this embedding the action of $\cuat$ on $\Sigma$ and
$\Sigma^*$ is:
\begin{center}
\begin{tabular}{lp{1cm}l}
$\bar{i} \cdot v_q=\left\{\begin{array}{ll} v_{q-1} & \textrm{if q
is even} \\ -v_{q+1} & \textrm{if q is odd}
\end{array}\right.$ & & $\bar{j}\cdot v_q=(-1)^{q+1}v_qi$ \vspace{3pt} \\
$v_q^* \cdot \bar{i}=\left\{\begin{array}{ll} -v_{q-1}^* & \textrm{if q is even} \\
v_{q+1}^* & \textrm{if q is odd}
\end{array}\right.$ & & $v_q^*\cdot \bar{j}=(-1)^{q+1}iv_q^*$ \\
\end{tabular}
\end{center}
These actions give rise to the following relations in $\Sigma^*
\otimes_{\cuat} \Sigma$. Let $p,q \in \{1,2,...,m\}$. Then:
$$\begin{array}{l}
v_{2p}^* \otimes v_{2q}=-v_{2p}^* \otimes \bar{i}\cdot
v_{2q-1}=-v_{2p}^*\cdot \bar{i} \otimes v_{2q-1}=v_{2p-1}^* \otimes v_{2q-1}, \\
v_{2p-1}^* \otimes v_{2q}=-v_{2p-1}^* \otimes \bar{i}\cdot
v_{2q}=v_{2p-1}^*\cdot \bar{i} \otimes v_{2q-1}=-v_{2p}^* \otimes
v_{2q-1}.
\end{array}$$
The set $\{v_{2p} \otimes v_l \vert p=1,...,m;l=1,...,n\}$ is a
basis of $\Sigma^* \otimes_{\cuat} \Sigma$ as a right $\Co$-vector
space. The left $\Co$-action on this coring is:
$$\begin{array}{ll}
i(v_{2p}^* \otimes v_{2q}) & =iv_{2p}^* \otimes
v_{2q}=(-1)^{2p+1}v_{2p}^*\cdot \bar{j} \otimes
v_{2q}=(-1)^{2p+1}v_{2p}^* \otimes \bar{j} \cdot v_{2q}
\\ & =v_{2p}^* \otimes v_{2q}i=(v_{2p}^* \otimes v_{2q})i, \vspace{3pt} \\
i(v_{2p}^* \otimes v_{2q-1}) & =iv_{2p}^* \otimes
v_{2q-1}=(-1)^{2p+1}v_{2p}^*\cdot \bar{j} \otimes
v_{2q-1}=(-1)^{2p+1}v_{2p}^* \otimes \bar{j} \cdot v_{2q-1}
\\ & =-v_{2p}^* \otimes v_{2q}i=-(v_{2p}^* \otimes v_{2q-1})i.
\end{array}$$
The comultiplication and counit of $\Sigma^* \otimes_{\cuat}
\Sigma$ is:
$$\begin{array}{ll}
\Delta(v_{2p} \otimes v_{2q}) & =\sum_{l=1}^m (v_{2p}^* \otimes
v_{2l}) \otimes_{\Co} (v_{2l}^* \otimes v_{2q})+ \sum_{l=1}^m
(v_{2p}^* \otimes v_{2l-1}) \otimes_{\Co} (v_{2l-1}^* \otimes
v_{2q}) \\
 & =\sum_{l=1}^m (v_{2p}^* \otimes v_{2l}) \otimes_{\Co} (v_{2l}^*
\otimes v_{2q})- \sum_{l=1}^m (v_{2p}^* \otimes v_{2l-1})
\otimes_{\Co} (v_{2l}^* \otimes v_{2q-1})\\
\epsilon(v_{2p} \otimes v_{2q}) & =\delta_{p,q} \vspace{3pt} \\
\Delta(v_{2p} \otimes v_{2q-1}) & =\sum_{l=1}^m (v_{2p}^* \otimes
v_{2l}) \otimes_{\Co} (v_{2l}^* \otimes v_{2q-1})+ \sum_{l=1}^m
(v_{2p}^* \otimes v_{2l-1}) \otimes_{\Co} (v_{2l-1}^* \otimes
v_{2q-1}) \\
 & =\sum_{l=1}^m (v_{2p}^* \otimes v_{2l}) \otimes_{\Co} (v_{2l}^*
\otimes v_{2q})+ \sum_{l=1}^m (v_{2p}^* \otimes v_{2l-1})
\otimes_{\Co} (v_{2l}^* \otimes v_{2q})\\
\epsilon(v_{2p} \otimes v_{2q-1}) & =0. \\
\end{array}$$
Let $\coring{T}=\Co c \oplus \Co s$ denote the trigonometric
coring and let $M^c(\Real,m)$ be the $\Real$-comatrix coalgebra of
order $m$. Then $\coring{T} \otimes_{\Real} M^c(\Real,m)$ becomes
a $\Co$-coring in the natural way \cite{Gomez/Louly:2003}. It may
be verified that the map from $\coring{T} \otimes_{\Real}
M^c(\Real,m)$ to $\Sigma^* \otimes_{\cuat} \Sigma$ defined by
$$c \otimes_{\Real} x_{pq} \mapsto v_{2p}^* \otimes v_{2q}, \qquad
s \otimes_{\Real} x_{pq} \mapsto v_{2p}^* \otimes v_{2q-1},$$ is
an isomorphism of corings.

\end{document}